\documentclass{birkmult}

\usepackage{amsfonts}
\usepackage{graphicx}
\usepackage{amsmath}
\usepackage{amssymb}

 \advance\hoffset by -0.0 truecm

\newtheorem{theorem}{Theorem}

\newtheorem{proposition}{Proposition}
\theoremstyle{definition}
\newtheorem{definition}{Definition}
\theoremstyle{remark}
\newtheorem{remark}{Remark}
\numberwithin{equation}{section}

\DeclareMathOperator{\spn}{span}

\DeclareMathOperator{\cotan}{cotan}

\begin{document}

\title[Sub-Riemannian geodesics]
{Sub-Riemannian geodesics on the 3-D sphere}

\author{Der-Chen Chang}

\address{Department of Mathematics, Georgetown University, Washington
D.C. 20057, USA}

\email{chang@georgetown.edu}

\thanks{}

\author{Irina Markina}
\address{Department of Mathematics,
University of Bergen, Johannes Brunsgate 12, Bergen 5008, Norway}

\email{irina.markina@uib.no}

\author{Alexander Vasil'ev}

\address{Department of Mathematics,
University of Bergen, Johannes Brunsgate 12, Bergen 5008, Norway}

\email{alexander.vasiliev@uib.no}

\thanks{The first author is partially supported by a research
grant from the United State Army Research Office and a Hong Kong RGC
competitive earmarked research grant $\#$600607. The second and the
third authors have been  supported by the grant of the Norwegian Research Council \#177355/V30, and by the European Science Foundation Research Networking Programme HCAA}




\dedicatory{to Bj\"orn Gustafsson on the occasion of his 60-th birthday}


\begin{abstract}
The unit sphere $\mathbb S^3$ can be identified with the unitary group $SU(2)$. Under this identification the unit sphere can be considered as a non-commutative Lie group. The commutation relations for the vector fields of the corresponding Lie algebra define a 2-step sub-Riemannian manifold. We study sub-Riemannian geodesics on this sub-Riemannian manifold making use of  the Hamiltonian formalism and solving the corresponding Hamiltonian system.

\end{abstract}

\maketitle

\section{Introduction}

The unit $3$-sphere centered on the origin is the set of $\mathbb
R^4$ defined by
$$\mathbb S^3=\{(x_1,x_2,x_3,x_4)\in\mathbb R^4:\ x_1^2+x_2^2+x_3^2+x_4^2=1)\}.$$
It is often convenient to regard $\mathbb R^4$ as the two complex dimensional space $\mathbb C^2$ or the space of
quaternions $\mathbb H$. The unit 3-sphere is then given by
$$\mathbb S^3=\{(z_1,z_2)\in\mathbb C^2:\ |z_1|^2+|z_2|^2=1)\}\quad\text{or}\quad
\mathbb S^3=\{q\in\mathbb H:\ |q|^2=1)\}.$$ The latter description represents the sphere $\mathbb S^3$ as a set of
unit quaternions and it can be considered as a group $Sp(1)$,
where the group operation is just a multiplication of quaternions. The group $Sp(1)$ is a
three-dimensional Lie group, isomorphic to $SU(2)$ by the isomorphism
$\mathbb C^2\ni (z_1,z_2)\to q\in\mathbb H$. The unitary group $SU(2)$ is the group of matrices
\begin{equation*}
\begin{pmatrix}
  z_1   &  z_2   \\
  -\bar z_2  & \bar z_1
\end{pmatrix}, \qquad z_1, z_2 \in \mathbb C,\quad |z_1|^2+|z_2|^2=1,\end{equation*} where the group law is given
by the multiplication of matrices. Let us identify $\mathbb R^3$ with
pure imaginary quaternions. The conjugation $qh\bar q$ of a pure
imaginary quaternion $h$ by a unit quaternion $q$ defines 
rotation in $\mathbb R^3$, and since $|qh\bar q|=|h|$, the map
$h\mapsto qh\bar q$ defines a two-to-one homomorphism $Sp(1)\to
SO(3)$. The Hopf map $\pi:\mathbb S^3\to\mathbb S^2$ can be
defined by $$\mathbb S^3\ni q\mapsto qi\bar q=\pi(q)\in\mathbb
S^2.$$ The Hopf map defines a principle circle bundle also known
as the Hopf bundle. Topologically $\mathbb S^3$ is a compact,
simply-connected, 3-dimensional manifold without boundary.

Even a small part of properties of the unit
3-sphere finds numerous applications in complex geometry,
topology, group theory, mathematical physics and others fields 
of mathematics. In the present paper we give a new look at the
unit 3-sphere, considering it as a sub-Riemannian manifold. The
sub-Riemannian structure comes naturally from the non-commutative
group structure of the sphere in a sense that two vector fields
span the smoothly varying distribution of the tangent bundle and
their commutator generates the missing direction. The
sub-Riemannian metric is defined as the restriction of the euclidean
inner product from $\mathbb R^4$ to the distribution. The present
paper devoted to the description of sub-Riemannian geodesics
on the sphere. The sub-Riemannian geodesics are defined as a
projection of the solution to the corresponding
Hamiltonian system onto the manifold. We give explicit formulas using 
different parametrizations and discuss the number of geodesics
starting from the unity of the group. While working on this paper the authors became
aware on the results in \cite{Hurtado} where the Lagrangian approach was developed and
the minimizers were found (in our terminology geodesics are solutions to a Hamiltonian system so
a minimizer is one of them).

\section{Left-invariant vector fields and the horizontal distribution}
In order to calculate left-invariant vector fields we use the definition of $\mathbb S^3$ as a set of unit
quaternions equipped with the following noncommutative multiplication ``$\circ$'': if $x=(x_1,x_2,x_3,x_4)$ and
$y=(y_1,y_2,y_3,y_4)$, then
\begin{eqnarray}\label{grlaw2}
x\circ y=(x_1,x_2,x_3,x_4)\circ(y_1,y_2,y_3,y_4) & = & \Big((x_1y_1-x_2y_2-x_3y_3-x_4y_4),
\nonumber \\ & \ \ & (x_2y_1+x_1y_2-x_4y_3+x_3y_4),
\nonumber
\\& \ & (x_3y_1+x_4y_2+x_1y_3-x_2y_4),\\ & \ & (x_4y_1-x_3y_2+x_2y_3+x_1y_4)\Big).\nonumber
\end{eqnarray}
The rule~\eqref{grlaw2} gives us the left translation $L_x(y)$ of an element $y=(y_1,y_2,y_3,y_4)$ by the element $x=(x_1,x_2,x_3,x_4)$. The left-invariant basis vector fields are defined as $X(x)=(L_x(y))_*X(0)$, where $X(0)$ are basis vectors at the unity of the group. The matrix corresponding to the tangent map $(L_x(y))_*$ calculated by~\eqref{grlaw2} becomes
$$(L_x(y))_* =\begin{pmatrix}
 x_1  & -x_2   & -x_3 & -x_4\\
x_2  & x_1   & -x_4 & x_3\\
x_3  & x_4  & x_1 & -x_2\\
x_4  & -x_3   &x_2 & x_1
\end{pmatrix}. $$
Calculating the action of $(L_x(y))_*$ in the basis of unit vectors of $\mathbb R^4$ we get four vector fields
\begin{eqnarray}\label{vf}
X_1(x) & = & x_1\partial_{x_1}+x_2\partial_{x_2}+x_3\partial_{x_3}+x_4\partial_{x_4},\nonumber \\
X_2(x) & = & -x_2\partial_{x_1}+x_1\partial_{x_2}+x_4\partial_{x_3}-x_3\partial_{x_4},\\
X_3(x) & = & -x_3\partial_{x_1}-x_4\partial_{x_2}+x_1\partial_{x_3}+x_2\partial_{x_4},\nonumber \\
X_4(x) & = & -x_4\partial_{x_1}+x_3\partial_{x_2}-x_2\partial_{x_3}+x_1\partial_{x_4}.\nonumber
\end{eqnarray} It is easy to see that the vector $X_1(x)$ is the unit normal to $\mathbb S^3$ at $x$ with respect to the usual inner product $\langle\cdot,\cdot\rangle$ in $\mathbb R^4$, hence, we denote  $X_1(x)$ by $N$. Moreover, $$\langle N, X_2(x)\rangle=\langle N, X_3(x)\rangle=\langle N, X_4(x)\rangle=0,\quad\text{and}\quad|X_k(x)|^2=\langle X_k(x), X_k(x)\rangle=1,$$ for $k=2,3,4$, and for any $x\in\mathbb S^3$. The matrix
$$\begin{pmatrix}
 -x_2  & x_1   & x_4 & -x_3\\
-x_3  & -x_4  & x_1 & x_2\\
-x_4  & x_3   &-x_2 & x_1
\end{pmatrix}$$ has rank 3, and we conclude that the vector fields $X_2(x)$, $X_3(x)$, $X_4(x)$ form an orthonormal basis  of the tangent space $T_x\mathbb S^3$ with respect to $\langle\cdot,\cdot\rangle$ at any point $x\in\mathbb S^3$. Let us denote the vector fields by $$X_3=X,\quad X_4=Y,\quad X_2=Z.$$

The vector fields possess the following commutation relations
$$[X,Y]=XY-YX=2Z,\quad [Z,X]=2Y,\quad [Y,Z]=2X.$$ Let $\mathcal{D} = \spn
\{ X, Y\}$ be the distribution generated by the vector fields $X$
and $Y$. Since $[X, Y] =  2Z \notin \mathcal{D} $, it follows that
$\mathcal{D}$ is not involutive. The distribution $\mathcal{D}$
will be called  {\it horizontal}. Any curve on the sphere with the
velocity vector contained in the distribution $\mathcal{D}$ will
be called a {\it horizontal curve}. Since $T_x\mathbb
S^3=\spn\{X,Y,Z=1/2[X,Y]\}$, the distribution is bracket
generating. We define the metric on the distribution $\mathcal D$
as the restriction of the metric $\langle\cdot,\cdot\rangle$ onto
$\mathcal D$, and the same notation $\langle\cdot,\cdot\rangle$ will be used. The manifold
$(\mathbb S^3,\mathcal D,\langle\cdot,\cdot\rangle)$ is a step two
sub-Riemannian manifold.

\begin{remark}
Notice that the choice of the horizontal distribution is not unique.
The relations $[Z,X]=2Y$ and $[Y,Z]=2X$ imply possible choices $\mathcal{D} =
\spn \{ X, Z\}$ or $\mathcal{D} = \spn \{ Y, Z\}$. The geometries
defined by different horizontal distributions are cyclically
symmetric, so we restrict our attention to the $\mathcal{D} = \spn
\{ X, Y\}$.
\end{remark}

We also can define the distribution as a kernel of the following one form
$$\omega = -x_2 dx_1 + x_1 dx_2 + x_4 dx_3 - x_3 dx_4$$ on $\mathbb{R}^4$.
One can easily check that
$$
\omega(X) =0,\quad \omega(Y)=0,\quad \omega (Z)=1 \not=0,\quad \omega (N)=0.$$
Hence, $\ker
\omega = \spn\{ X, Y, N \}$, and the horizontal distribution can be
written as
$$\mathbb{S}^3 \ni x \rightarrow \mathcal{D}_x = \ker \omega \cap T_x \mathbb{S}^3. $$

Let $\gamma(s)=(x_1(s),x_2(s),x_3(s),x_4(s))$ be a curve on $\mathcal S^3$. Then the velocity vector, written in the left-invariant basis, is $$\dot\gamma(s)=a(s)X(\gamma(s))+b(s)Y(\gamma(s))+c(s)Z(\gamma(s)),$$ where
\begin{eqnarray}\label{coord}
a & = & \langle\dot\gamma, X\rangle=-x_3\dot x_1-x_4\dot x_2+x_1\dot x_3+x_2\dot x_4,\nonumber \\
b & = & \langle\dot\gamma, Y\rangle=-x_4\dot x_1+x_3\dot x_2-x_2\dot x_3+x_1\dot x_4, \\
c & = & \langle\dot\gamma, Y\rangle=-x_2\dot x_1 +x_1\dot x_2+x_4\dot x_3-x_3\dot x_4.\nonumber
\end{eqnarray}
The following proposition holds.

\begin{proposition}
\label{prop1} Let $\gamma (s) = (x_1(s), x_2(s), y_1(s),
y_2(s))$ be a curve on $\mathbb{S}^3$. The curve $\gamma$ is
horizontal, if and only if,
\begin{equation}\label{hc}
c=\langle \dot\gamma, Z \rangle  = \langle\dot\gamma, X\rangle=-x_2\dot x_1 +x_1\dot x_2+x_4\dot x_3-x_3\dot x_4=0.
\end{equation}
\end{proposition}

The manifold $\mathbb S^3$ is connected and it satisfies the
bracket generating condition. By the Chow theorem~\cite{Chow},
there exists piecewise $C^1$ horizontal curves connecting two
arbitrary points on $\mathbb S^3$. In fact, smooth horizontal
curves connecting two arbitrary points on $\mathbb S^3$ were
constructed in~\cite{CChM}.

\begin{proposition}
The horizontality property is invariant under the left translation.
\end{proposition}

\begin{proof}
It can be shown that~\eqref{coord} does not change under the left translation. This implies the conclusion of the
proposition.
\end{proof}

\section{Hamiltonian system}

Ones we have a system of curves, in our case the system of horizontal curves, we can define the length as in
the Riemannian geometry. Let $\gamma:[0,1]\to \mathbb S^3$ be a horizontal curve such that $\gamma(0)=x$,
$\gamma(1)=y$, then the length $l(\gamma)$ of $\gamma$ is defined as the following
\begin{equation}\label{length}
l(\gamma)=\int_0^1\langle\dot\gamma, \dot\gamma\rangle^{1/2}\,dt=\int_0^1\big(a^2(t)+b^2(t)\big)^{1/2}\,dt.
\end{equation} Now we are able to define the distance between the points $x$ and $y$ by minimizing the
integral~\eqref{length} or the corresponding energy integral
$\int_0^1\big(a^2(t)+b^2(t)\big)\,dt$ under the non-holonomic
constraint~\eqref{hc}. This is a Lagrangian approach. The Lagrangian
formalism was applied to study the sub-Riemannian geometry of
$\mathbb S^3$ in~\cite{CChM, Hurtado}. In the Riemannian geometry the
minimizing curve locally coinsides with the geodesic, but it is not
the case for the sub-Riemannian manifolds. Interesting examples
and discussions can be found, for instance
in~\cite{LS,Mon1,Mon2,Mon,Str}. Given the sub-Riemannian metric we
can form a Hamiltonian function defined on the cotangent bundle of
$\mathbb S^3$. The geodesics in the sub-Riemannian manifolds are
defined as a projection of the solution to the
corresponding Hamiltonian system onto the manifold. It is a good generalization of
the Riemannian case in the following sense. The Riemannian
geodesics (that are defined as curves with vanishing acceleration) can be lifted
to the solutions of the Hamilton system on the cotangent bundle.

In the present paper we are interested in the construction of sub-Riemannian geodesics
on $(\mathbb S^3,\mathcal D,\langle\cdot,\cdot\rangle)$. Let us
write the left-invariant vector fields $X,Y,Z$, using the matrices
\begin{equation}
\label{eq:matrix}
I_1=\left[\array{rrrr}0 & 0 & -1 & 0
\\ 0 & 0 & 0 & -1
\\
1 & 0 & 0 & 0
\\ 0 & 1 & 0 & 0\endarray\right],\quad
I_2=\left[\array{rrrr}0 & 0 & 0 & -1
\\ 0 & 0 & 1 & 0
\\
0 & -1 & 0 & 0
\\ 1 & 0 & 0 & 0\endarray\right],\quad
I_3=\left[\array{rrrr}0 & -1 & 0 & 0
\\ 1 & 0 & 0 & 0
\\
0 & 0 & 0 & 1
\\ 0 & 0 & -1 & 0\endarray\right].
\end{equation} Then $$X=\langle I_1x,\nabla x\rangle,\quad Y=\langle I_2x,\nabla x\rangle,
\quad Z=\langle I_3x,\nabla x\rangle.$$
The Hamiltonian function is defined as
$$H=\frac{1}{2}(X^2+Y^2)=\frac{1}{2}\Big(\langle I_1x,\xi\rangle^2
+\langle I_2x,\xi\rangle^2\Big),$$ where $\xi=\nabla x$. Then the Hamiltonian system follows as
\begin{equation}\label{hs}
\begin{split}
&\dot x=\frac{\partial H}{\partial \xi}\ \ \Rightarrow \ \ \dot
x=\langle I_1x,\xi\rangle\cdot (I_1x)+\langle I_2x,\xi\rangle\cdot (I_2x)\\
&\dot \xi=-\frac{\partial H}{\partial x}\ \ \Rightarrow \ \
\dot\xi=\langle I_1x,\xi\rangle\cdot (I_1\xi)+\langle
I_2x,\xi\rangle\cdot (I_2\xi).
\end{split}\end{equation}
As it was mentioned, a geodesic is the projection of a solution to the
Hamiltonian system onto the $x$-space. We obtain the following
properties.
\begin{itemize}
 \item [1.]{
Since $\langle I_1x,x\rangle=\langle I_2x,
x\rangle=\langle I_3x,x\rangle=0$, multiplying the first equation
of~\eqref{hs} by $x$ we get
$$\langle\dot x,x\rangle=0\ \ \Rightarrow \ \ |x|^2=const.$$
We conclude that {\it any solution to the Hamiltonian system belongs to the
sphere.} Taking the constant equal to $1$ we get geodesics on
$\mathbb S^3$.}
\item[2.]{Multiplying the first equation of~\eqref{hs} by
$I_3x$, we get
\begin{equation}\label{hc1}
\langle\dot x,I_3x\rangle=0,
\end{equation}
by the rule of multiplication for $I_1$, $I_2$, and $I_3$. The
reader easily recognizes the horizontality condition $\langle\dot
x,Z\rangle=0$ in~\eqref{hc1}. It means that {\it any solution to
the Hamiltonian system is a horizontal curve}}.
\item[3.]{Multiplying the first equation of~\eqref{hs} by $I_1x$, and then by $I_2x$,
we get
$$\langle \xi,I_1x\rangle=\langle\dot x,I_1x\rangle,
\qquad \langle \xi,xI_2\rangle=\langle\dot x,xI_2\rangle.$$ From
the other side, we know that $\langle\dot x,I_1x\rangle=a$ and
$\langle\dot x,xI_2\rangle=b$. The Hamiltonian function can be
written in the form
$$H=\frac{1}{2}\Big(\langle I_1x,\xi\rangle^2+\langle I_2x,\xi\rangle^2\Big)
=\frac{1}{2}\Big(\langle I_1x,\dot x\rangle^2+\langle I_2x,\dot
x\rangle^2\Big)=\frac{1}{2}\Big(a^2+b^2\Big).
$$
Thus, {\it the Hamiltonian function gives the kinetic energy $H=\frac{|\dot
q|^2}{2}$ and it is a constant along the geodesics.}}
\item[4.]
{If we multiply the first equation of~\eqref{hs}
by $\dot x$, then we get
\begin{equation*}
|\dot x|^2=\langle I_1x,\xi\rangle^2+\langle I_2x,\xi\rangle^2
=\langle I_1x,\dot x\rangle^2+\langle I_2x,\dot
x\rangle^2=a^2+b^2=2H.
\end{equation*}  Therefore \begin{equation}\label{eq1}|\dot
x|^2=a^2+b^2.\end{equation}}
\end{itemize}

\section{Velocity vector with constant coordinates}

We know that the length of the velocity vector is
constant along geodesics. Let us start from the simplest case, when the coordinates of
the velocity vector are constant. Suppose that $\dot a=\dot b=0$. The first line of system~\eqref{hs}
can be written as 
\begin{eqnarray}\label{hs1}
\dot x_1 & =  -ax_3-bx_4 \qquad \qquad \dot x_3 & =  +ax_1-bx_2  \\
\dot x_2 & =  -ax_4+bx_3 \qquad \qquad \dot x_4 & =  +ax_2+bx_1 \nonumber.
\end{eqnarray}

Differentiation of system~\eqref{hs1} yields
\begin{eqnarray}\label{hs2}
\ddot x_1 & =  -a\dot x_3-b\dot x_4 \qquad \qquad \ddot x_3 & =  +a\dot x_1-b\dot x_2 \\
\ddot x_2 & =  -a\dot x_4+b\dot x_3  \qquad\qquad \ddot x_4 & =  +a\dot x_2+b\dot x_1 \nonumber.
\end{eqnarray} We substitute the first derivatives from~\eqref{hs1} in~\eqref{hs2}, and get
\begin{equation}\label{hs3}
\ddot x_k=-r^2 x_k,\quad r^2=a^2+b^2,\quad k=1,2,3,4.
\end{equation}

\begin{theorem}\label{th1}
The set of geodesics with constant velocity coordinates form a
unit sphere $\mathbb S^2$ in $\mathbb R^3$
\end{theorem}

\begin{proof}
We are looking for horizontal geodesics parametrized by the arc
length and starting from the point $x(0)=x_0$. So, we set $r=1$
and $a=\cos\psi$, $b=\sin\psi$, where $\psi$ is a constant from
$[0,2\pi)$. Solving the equation~\eqref{hs3} we get the general
solution $x(s)=A\cos s+B\sin s$. We conclude that $A=x_0$  from
the initial data. To find $B$ let us substitute the general
solution in equations~\eqref{hs1} and get $B=(aI_1+bI_2)x_0$.
Thus, the horizontal geodesics with constant horizontal
coordinates are $$x(s)=x_0\cos s + (\cos\psi I_1+\sin \psi
I_2)x_0\sin s$$ Since the geodesics are invariant under the left
translation it is sufficient to describe the situation at the
unity element, e.g., $x_0=(1,0,0,0)$ of $\mathbb S^3$. In this
case the geodesics are
\begin{eqnarray}\label{constgeod}
x_1 & =  \cos s,\qquad\qquad x_3 & =  \cos\psi\sin s, \\
x_2 & =  0, \qquad\qquad\ \ \ \  x_4 & =  \sin\psi\sin s. \nonumber
\end{eqnarray}
We see that the set of geodesics with constant velocity
coordinates form the unit sphere $\mathbb S^2$ in $\mathbb
R^3=\{(x_1,0,x_3,x_4)\}$. The parameter $\psi\in [0,2\pi)$
corresponds to the initial velocity.
\end{proof}
The sphere~\eqref{constgeod} is a direct analogue of the
horizontal plane in the Heisenberg group $\mathbb H^1$ at the unity. 
We remark that this result was obtained independently in \cite{Hurtado}.

Let us
calculate the analogue of the vertical axis in $\mathbb H^1$. We wish
to find an integral curve for the vector field $Z$. In other words,
we solve the system
\begin{eqnarray}\label{vertical}
a & = & \langle\dot\gamma, X\rangle=-x_3\dot x_1-x_4\dot x_2+x_1\dot x_3+x_2\dot x_4=0,\nonumber \\
b & = & \langle\dot\gamma, Y\rangle=-x_4\dot x_1+x_3\dot x_2-x_2\dot x_3+x_1\dot x_4=0, \\
c & = & \langle\dot\gamma, Z\rangle=-x_2\dot x_1 +x_1\dot x_2+x_4\dot x_3-x_3\dot x_4=1,\nonumber \\
n & = & \langle\dot\gamma, N\rangle=+x_1\dot x_1 +x_2\dot
x_2+x_3\dot x_3+x_4\dot x_4=0.\nonumber
\end{eqnarray}
The determinant of the system is $1$ and it is redused to
\begin{eqnarray*}
\dot x_1 & =  -x_2,\qquad\qquad \dot x_3 & =  +x_4,\\
\dot x_2 & =  +x_1,\qquad\qquad \dot x_4 & =  -x_3.
\end{eqnarray*} Differentiating again, we get the equation $\ddot x=-x$.
The initial point is $x(0)=x_0$. System~\eqref{vertical} gives the value of the initial velocity
$\dot x(0)=I_3x_0$. Taking into account this initial data
we get the equation of the {\it vertical line} as $$x(s)=x_0\cos s +
I_3x_0\sin s.$$ In particular, at the point $(1,0,0,0)$ the
equation of the vertical line is
\begin{equation}\label{vcircl}
x_1  =  \cos s,\quad x_2  =  \sin s,\quad
x_3  =  0,\quad
x_4  =  0,\qquad s\in[0,2\pi].
\end{equation}

\section{Velocity vector with non-constant coordinates}

\noindent{\it Cartesian coordinates}

\medskip
Fix the initial point $x^{(0)}=(1,0,0,0)$. It is convenient to introduce complex coordinates $z=x_1+ix_2$, $w=x_3+ix_4$, $\varphi=\xi_1+i\xi_2$, and $\psi=\xi_3+i\xi_4$. Hence, the Hamiltonian admits the form $2H=|\bar w\varphi -z\bar\psi|^2$. The corresponding Hamiltonian system becomes
\begin{equation*}
\begin{array}{lllll}
\dot z & = &  w(\bar{w} \varphi- z\bar{\psi}),\quad & &z(0)=1, \\
\dot w & = & -z(w\bar{\varphi}-\bar{z}\psi), \quad & &w(0)=0,\\
\dot {\bar{\varphi}} & = & \bar{\psi}(w\bar{\varphi}-\bar{z}\psi), \quad & &\bar{\varphi}(0)=A-iB,\\
\dot {\bar{\psi}} & = & -\bar{\varphi}(\bar{w} \varphi-z\bar{\psi}), \quad & &\bar{\psi}(0)=C-iD.
\end{array}
\end{equation*}
Here the constants $B,C$, and $D$ have the following dynamical meaning: $\dot w(0)=C+iD$,
and $B=-i \ddot w(0)/2\dot w(0)$. So $C,D$ is the velocity and $B/\sqrt{C^2+D^2}$ is the curvature of a geodesic at the initial point. This complex Hamiltonian system has the
first integrals
\[
\begin{array}{rcl}
z\psi -w \varphi  & = &  C+iD, \\
z\bar{\varphi}+w \bar{\psi} & = & A-iB,\\
\end{array}
\]
and we have $|z|^2+|w|^2=1$ and $2H=C^2+D^2=1$ as an additional
normalization. The latter means that we parametrize geodesics by the natural parameter. Therefore,
\[
\begin{array}{rcl}
\varphi  & = & z(A+iB)-\bar w ( C+iD), \\
\psi  & = & \bar z (C+iD) +w (A+iB).\\
\end{array}
\]
Let us introduce an auxiliary function $p=\bar w/{z}$. Then substituting $\varphi$ and $\psi$ in the
Hamiltonian system we get the equation for $p$ as
\[
\dot{p}=(C+iD)p^2-2iBp+(C-iD),\quad p(0)=0.
\]
The solution is
\[
p(s)=\frac{(C-iD)\sin(s\sqrt{1+B^2})}{\sqrt{1+B^2}\cos(s\sqrt{1+B^2})+iB \sin(s\sqrt{1+B^2})}.
\]
Taking into account that $\dot z \bar z= -w \dot{\bar w}$, we get the solution
\begin{equation}\label{geodz}
z(s)=\left(\cos(s\sqrt{1+B^2})+i\frac{B}{\sqrt{1+B^2}}\sin(s\sqrt{1+B^2})\right)e^{-iBs},
\end{equation}
and
\begin{equation}\label{geodw}
w(s)=\frac{C+iD}{\sqrt{1+B^2}} \sin(s\sqrt{1+B^2})e^{iBs}.
\end{equation}
If $B=0$ we get the solutions with constant horizontal velocity coordinates $$z(s)=\cos s,
\qquad w(s)=(\dot x_3(0)+i\dot x_4(0))\sin s$$ from the previous section.

\begin{theorem}\label{numbergeod1}
Let $A$ be a point of the vertical line, i.~e. $A=(\cos \omega,\sin \omega, 0,0)$, $\omega\in[0,2\pi)$, then there are countably
many geometrically different geodesics $\gamma_n$ connecting $O=(1,0,0,0)$ with $A$. They have the following parametric equations
\begin{eqnarray}\label{geod1}
 z_n(s) & = & \Big(\cos\big(s\frac{\pi n}{\sqrt{\pi^2 n^2-\omega^2}}\big)-
 i\frac{\omega}{\pi n}\sin\big(s\frac{\pi n}{\sqrt{\pi^2 n^2-\omega^2}}\big)\Big)e^{\frac{is\omega}{\sqrt{\pi^2n^2-\omega^2}}},\\
w_n(s) & = & \big(\dot x_3(0)+i\dot x_4(0)\big)\frac{\sqrt{\pi^2n^2-\omega^2}}{\pi n}
\sin\big(s\frac{\pi n}{\sqrt{\pi^2 n^2-\omega^2}}\big)e^{\frac{-is\omega}{\sqrt{\pi^2n^2-\omega^2}}},\nonumber
\end{eqnarray}
$n\in\mathbb Z\setminus\{0,\pm 1\}$, $s\in[0,s_n]$, where $l_n=\frac{1}{\sqrt{2}}s_n=\frac{1}{\sqrt{2}}\sqrt{\pi^2 n^2-\omega^2}$ is the length of the geodesic
$\gamma_n$.
\end{theorem}
\begin{proof}
Since we use the condition $2H=|\dot z|^2+|\dot w|^2=1$ we
conclude that the geodesics are parametrized by the arc length, and the
length of a geodesic at the value of the parameter $s=l \sqrt{2}$ is equal to
$l$. If the point $A=(z(s),w(s))$ belongs to the vertical line starting at
$O=(1,0,0,0)$, then $|z(s)|=1$ and $|w(s)|=0$ provided that
$-Bs=\omega$. It implies $$\cos^2
(s\sqrt{1+B^2})+\frac{B^2}{1+B^2}\sin^2(s\sqrt{1+B^2})=1,\qquad
\sin(s\sqrt{1+B^2})=0,\qquad -Bs=\omega.$$ These equations are satisfied when
$$s_n=\sqrt{\pi^2 n^2-\omega^2},\quad B_n=-\frac{\omega}{\sqrt{\pi^2
n^2-\omega^2}},\quad n\in\mathbb Z\setminus\{0\}.$$ We conclude,
that for $n(s)\in\mathbb Z\setminus\{0\}$ there is a
constant $B_n(s)=-\frac{\omega}{\sqrt{\pi^2 n^2-\omega^2}}$, such that the
corresponding geodesic $\gamma_n(s)$, $s\in[0,s_n]$, satisfying 
equation~\eqref{geod1} joins the points $O$ and $A$ and the
length of the geodesic is equal to $s_n=\sqrt{\pi^2 n^2-\omega^2}$.
\end{proof}

\begin{remark}
In the formulation of the theorem the words `geometrically different' mean that due to
the rotation of the argument of $C+iD$ in $w(s)$, there exist uncountably many geodesics.
\end{remark}

So far we have had a clear picture of trivial geodesics whose velocity
has constant coordinates. They are essentially unique (up to
periodicity). The situation with geodesics joining the point
$(1,0,0,0)$ with the points of the vertical line $A$ has been
described in the preceding theorem. Let us consider the general
position of points on $\mathbb S^3$.

\begin{theorem}\label{numbergeod2}
Given an arbitrary point $(z_1,w_1)\in \mathbb S^3$ which neither belongs to the vertical line $A$ nor to
the horizontal sphere $\mathbb S^2$, there is a finite number of geometrically different geodesics
joining the initial point $(z_0,w_0)\in \mathbb S^3$ with $(z_1,w_1)$, $z_0=1$, $w_0=0$.
\end{theorem}
\begin{proof}
Let us denote
\[
w_1=\rho e^{i\varphi},\quad z_1=re^{i\alpha},\quad C+iD=e^{i\theta}.
\]
Then from \eqref{geodz} and \eqref{geodw} we have that
\begin{equation}\label{ugol}
r^2=1-\frac{1}{1+B^2}\sin^2(s\sqrt{1+B^2}),\quad \mbox{and\  \ }\varphi=Bs+\theta,
\end{equation} where $s$ is the value of the length arc parameter when the point $(z_1,w_1)$ is reached. We want to exclude the parameter $s$ and rewrite the equations~\eqref{geodz} and~\eqref{geodw} in terms of the parameter $B$ and the given dates $r,\rho\neq 0,\alpha\neq 0$, and $\varphi$.
We suppose for the moment that the angles $s\sqrt{1+B^2}$ and $sB$ are from the first quadrant.
Other cases are treated similarly. Then we have
\[
z=(\sqrt{1-(1+B^2) \rho^2}+iB\rho)e^{i(\theta-\varphi)},
\]
and
\[
\theta=\theta(B)=\alpha+\varphi-\arctan\frac{B\rho}{\sqrt{1-(1+B^2)\rho^2}}.
\]
The first expression in \eqref{ugol} leads to the value of the length parameter $s$ at $(z_1,w_1)$
\[
s=\frac{1}{\sqrt{1+B^2}}\arcsin (\rho\sqrt{1+B^2}),
\]
and the second to
\[
\varphi=\theta+\frac{B}{\sqrt{1+B^2}}\arcsin (\rho\sqrt{1+B^2}).
\]
Substituting $\theta(B)$ in the latter equation we obtain
\begin{equation}\label{equa}
\sin\left((\alpha-\arctan\Big(\frac{B\rho}{r^2-B\rho^2}\Big))\sqrt{1+\frac{1}{B^2}}\right)=\rho\sqrt{1+B^2},
\end{equation}
as an equation for the parameter $B$. Observe that $\varphi-\theta(B)=\alpha-\arctan\Big(\frac{B\rho}{r^2-B\rho^2}\Big)$ is a bounded function and
$\lim_{B\to 0}\theta(B)\neq 0$. Indeed, if the latter limit were vanishing, then the value of given $\varphi$
would be zero and the solution of the problem would be only $B=0$ which is the trivial case excluded
from the theorem. So the left-hand side of equation \eqref{equa} is a  function of $B$ which is bounded by 1 in
absolute value and fast oscillating about the point $B=0$. Observe, that $\alpha=0$ corresponds to the horizontal sphere which also was excluded from the theorem. The right-hand side of \eqref{equa}
is an even function increasing for $B>0$, see Figure~\eqref{fig2}. Therefore, there exists a countable number of non-vanishing
different solutions $\{B_n\}$ of the equation \eqref{equa} within the interval $|B|\leq \sqrt{\frac{1}{\rho^2}-1}=\frac{|z|}{|w|}$
with a limit point at the origin. 

\begin{figure}[ht]
\centering \scalebox{0.3}{\includegraphics{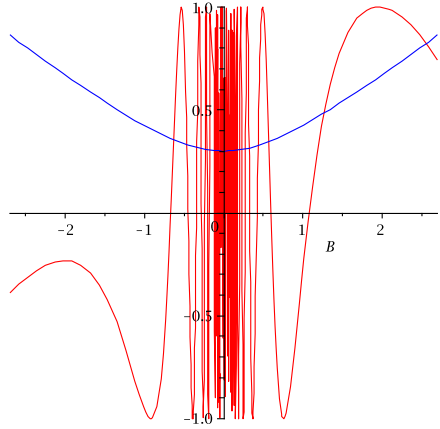}}
\caption[]{Solutions to the equation~\eqref{equa}  }\label{fig2}
\end{figure}

However, in order do define the parameters $B$,  we need
to solve the equations (\ref{ugol}), (\ref{equa}), and not all $B_n$ satisfy all three equations. Let us consider positive $B_n$.
We calculate the argument of $z$ as
\[
\alpha=-B_ns+\arctan\Bigg[\frac{B_n}{\sqrt{B_n^2+1}}\tan\left(s\sqrt{B_n^2+1}\right)\Bigg]
\]
\[
=-B_ns+\arctan\Bigg[\frac{B_n\rho}{\sqrt{1-(1+B_n^2)\rho^2}}\Bigg]
\]
\[
<-B_n s+\frac{B_n\rho}{\sqrt{1-(1+B_n^2)\rho^2}}.
\]

On the other hand, we have
\[
1=\frac{\arcsin(\rho\sqrt{1+B_n^2})}{s\sqrt{1+B_n^2}}>\frac{\rho}{s}.
\]
Observe that due the remark before this theorem, $\alpha>0$ and $0<\rho<1$.
Therefore, we deduce the inequality
\[
\alpha<B_n\rho\frac{1-\sqrt{1-\rho^2(1+B_n^2)}}{\sqrt{1-\rho^2(1+B_n^2)}},
\]
or
\begin{equation}\label{equ}
B_n\rho>\alpha \frac{\sqrt{1-\rho^2(1+B_n^2)}}{1-\sqrt{1-\rho^2(1+B_n^2)}}.
\end{equation}
The right-hand side of the inequality \eqref{equ} decreases with respect to $B_n>0$.

Set $\varepsilon=\frac{1+\rho^2}{2}$. If $\varepsilon <\rho^2(1+B_n^2)<1$, then immediately
we have the inequality
$B_n^2>\frac{1}{2}(\frac{1}{\rho^2}-1)>0$. If $0 <\rho^2(1+B_n^2)\leq \varepsilon$, then the inequality
\eqref{equ} implies that
\[
B_n>\alpha\frac{\sqrt{1-\varepsilon}}{\rho(1-\sqrt{1-\varepsilon})}=\alpha\frac{\sqrt{1-\rho^2}}{\rho(\sqrt{2}-\sqrt{1-\rho^2})}>0.
\]
Finally, we obtain
\[
B_n>\min\left\{\alpha\frac{\sqrt{1-\rho^2}}{\rho(\sqrt{2}-\sqrt{1-\rho^2})}, \sqrt{\frac{1}{2}(\frac{1}{\rho^2}-1)}\right\}\equiv b(\xi_1,\rho)>0.
\]
This proves that all positive solutions to the equation \eqref{equa} must belong to the interval $(b(\xi_1,\rho), \sqrt{\frac{1}{\rho^2}-1})$, hence there are only finite number of such $B_n$. The same arguments are applied for negative values of $B_n$. 

Let us discuss the limiting cases. If $\rho\to 0$ then the endpoint aproaches the vertical line. In this case the graph of the right hand side function in~\eqref{equa} approaches the horizontal axis and the range of $|B|$  increases. If $\alpha\to 0$ then the end point aproaches the horizontal sphere $\mathbb S^2$, the number of geodesics is finite for any value $\alpha\neq 0$, but decreases. 
\end{proof}

This theorem reveals similarity of sub-Riemannian geodesics on the sphere with those for the Heisenberg group. The number of geodesics joining the origin with a point neither from the vertical axis nor from the horizontal plane is finite and approaching the vertical line becomes infinite.

\medskip

\noindent{\it Hyperspherical coordinates}

\medskip

Let us use the hyperspherical coordinates to find  geodesics
 with non-constant velocity coordinates.
\begin{eqnarray}\label{hspc}
x_1+ix_2 & = & e^{i\xi_1}\cos\eta,\\
x_3+ix_4 & = & e^{i\xi_2}\sin\eta,\qquad \eta\in[0,\pi/2),\quad\xi_1,\xi_2\in[-\pi,\pi)\nonumber
\end{eqnarray}

The horizontal coordinates are written as
\begin{eqnarray*}
a & = & \dot\eta\cos(\xi_1-\xi_2)+(\dot\xi_1+\dot\xi_2)\sin(\xi_1-\xi_2)\frac{\sin 2\eta}{2},\\
b & = & -\dot\eta\sin(\xi_1-\xi_2)+(\dot\xi_1+\dot\xi_2)\cos(\xi_1-\xi_2)\frac{\sin 2\eta}{2},\\
c & = & \dot\xi_1\cos^2\eta-\dot\xi_2\sin^2\eta.
\end{eqnarray*}

The horizontality condition in hyperspherical coordinates becomes
$$\dot\xi_1\cos^2\eta-\dot\xi_2\sin^2\eta=0.$$

The horizontal sphere~\eqref{constgeod} is obtained from the
parametrization~\eqref{hspc}, if we set $\xi_1=0$, $\xi_2=\psi$,
$\eta=s$. We get
$$a^2+b^2=1=\dot\eta^2\quad\Longrightarrow\quad a=\cos\psi,\quad
b=\sin\psi.$$ The vertical line is obtained from the
parametrization~\eqref{hspc} setting $\eta=0$, $\xi_1=s$.

Writing the vector fields $N,Z,X,Y$ in the hyperspherical
coordinates we get $$N=-2\cotan2\eta\partial_{\eta},\quad
Z=\partial_{\xi_1}-\partial_{\xi_2},$$
$$X=\sin(\xi_1-\xi_2)\tan\eta\partial_{\xi_1}+\sin(\xi_1-\xi_2)\cotan\eta\partial_{\xi_2}
+2\cos(\xi_1-\xi_2)\partial_{\eta},$$
$$Y=\cos(\xi_1-\xi_2)\tan\eta\partial_{\xi_1}
+\cos(\xi_1-\xi_2)\cotan\eta\partial_{\xi_2}-2\sin(\xi_1-\xi_2)\partial_{\eta}.$$
In this parametrization the similarity with the Heisenberg group
can be shown. The commutator of two horizontal vector fields $X,Y$ gives the constant vector field $Z$ which is orthogonal to the
horizontal vector fields at each point of the manifold. In
hyperspherical coordinates it is easy to see that the form
$\omega=\cos^2\eta d\xi_1-\sin^2\eta d\xi_2$, that defines the
horizontal distribution is contact because 
$$\omega\wedge d\omega=\sin (2\eta)\,d\eta\wedge d\xi_1\wedge
d\xi_2=2dV,$$ where $dV$ is the volume form. The sub-Laplacian is
defined as
$$\frac{1}{2}(X^2+Y^2)=\frac{1}{2}(\tan^2\eta\partial^2_{\xi_1}+\cotan^2\eta\partial^2_{\xi_2}
+4\partial^2_{\eta}+2\partial_{\xi_1}\partial_{\xi_2}).$$
The Hamiltonian becomes
$$H(\xi_1,\xi_2,\eta,\psi_1,\psi_2,\theta)=\frac{1}{2}(\tan^2\eta\psi^2_1+\cotan^2\eta\psi^2_2
+4\theta^2+2\psi_1\psi_2),$$
and the corresponding Hamiltonian system is given as 
\begin{eqnarray*}
\dot\xi_1 & = & \frac{\partial H}{\partial \psi_1}=\psi_1\tan^2\eta+\psi_2\\
\dot\xi_2 & = & \frac{\partial H}{\partial \psi_2}=\psi_2\cotan^2\eta+\psi_1\\
\dot\eta & = & \frac{\partial H}{\partial \theta}=4\theta\\
\dot\psi_1 & = & -\frac{\partial H}{\partial \xi_1}=0\\
\dot\psi_2 & = & -\frac{\partial H}{\partial \xi_2}=0\\
\dot\theta & = & -\frac{\partial H}{\partial \eta}=-\psi_1^2\frac{\tan\eta}{\cos^2\eta}
+\psi_2^2\frac{\cotan\eta}{\sin^2\eta}\\
\end{eqnarray*} Let us solve this Hamiltonian system for the following initial data: $\eta(0)=0$, $\xi_1(0)=0$,
$\xi_2(0)=0$, $\psi_1(0)=\psi_1$, $\psi_2(0)=\psi_2$, $\theta(0)=\frac{\dot\eta(0)}{4}=\theta_0$.

We see that $\psi_1$ and $\psi_2$ are constant. The horizontality condition at $(0,0,0)$ gives $\dot\xi(0)=0$ and the first equation of Hamiltonian system implies that $\psi_2=0$. If $\psi_1=0$, then
$\dot\xi_1=\dot\xi_2=0$, $\dot\eta=4\theta_0$ and we get the
variety of trivial geodesics~\eqref{constgeod} up to the
reparametrization $s\mapsto 4\theta_0 s$. To find other geodesics
we suppose that $\psi_1\neq 0$ and $\dot\eta(0)>0$. This condition ensures us that the trajectory  starting at the point $(0,0,0,)$ remains in the domain of  parametrizaition locally in time. From the third and from the last equations of the Hamiltonian system we have
$$\ddot\eta=-4\psi_1^2\frac{\sin\eta}{\cos^3\eta}\qquad\Rightarrow\qquad\dot\eta\,d\dot\eta=-4\psi_1^2\frac{\sin\eta}{\cos^3\eta}d\eta
\qquad\Rightarrow$$
$$\dot\eta^2=C-4\frac{\psi_1^2}{\cos^2\eta},$$
We observe that
$C=\dot\eta^2(0)+4\psi_1^2>0$.

Continue to solve the Hamiltonian system finding $\eta(s)$.
$$\frac{\cos\eta\,d\eta}{\sqrt{C\cos^2\eta-4\psi_1^2}}=ds.$$
Denote by $\sin\eta=p$. Then,
\begin{equation}\label{eq10}
\frac{dp}{\sqrt{-Cp^2+C-4\psi_1^2}}=ds.
\end{equation}
Integrating~\eqref{eq10} from $0$ to $s$ we get
$$s+K=\frac{1}{\sqrt{|C|}}\arcsin\Big(\sqrt{\frac{C}{C-4\psi_1^2}}\sin\eta(s)\Big),$$
where $K$ is found setting $s=0$ as
$K=\frac{1}{\sqrt{|C|}}\arcsin 0=0$. We calculate
\begin{equation}\label{eta}
\sin^2\eta(s)=\frac{C-4\psi^2_1}{C}\sin^2(\sqrt C s).
\end{equation} From the Hamiltonian system we find
\begin{equation}\label{xi2}\xi_2(s)=\psi_1s,
\end{equation}
and $$\dot\xi_1=\psi_1\frac{\sin^2\eta(s)}{1-\sin^2\eta(s)}
=\psi_1\frac{\sin^2(\sqrt C s)}{a+\sin^2(\sqrt C s)},\quad
a=\frac{C}{C-4\psi_1^2}.$$ It gives
\begin{equation}\label{xi1}\xi_1(s)=-\psi_1s+\frac{\psi_1}{2|\psi_1|}\arctan\Big[\frac{2|\psi_1|}{\sqrt
C}\tan(\sqrt C s)\Big].
\end{equation}

Let us suppose for the moment that the geodesics are parametrized on the interval
$[0,1]$. If the initial point and the finite point are on the
vertical line: $\eta(0)=\eta(1)=0$, then
$$0=\sin^2\eta(1)=\frac{C-4\psi^2_1}{C}\sin^2(\sqrt C)\quad\Rightarrow\quad
C=\pi^2n^2.$$ Since the value of $\xi_2$ on the vertical line is arbitrary, the values of $C$ and $\xi_1(1)$ give us the value
of $\xi_2(1)$. Setting $C=\pi^2n^2$ in the equation for
$\xi_1(1)$, we find $\psi_1=-\xi_1(1)$. Then
$$\xi_2(1)=-\xi_1(1).$$ The finite point on the vertical
line corresponds to the value of $\xi_1(1), \eta(1)=0$, and
$\xi_2=-\xi_1(1)$.

We also note that the square of the velocity
$|v|^2=\dot\eta^2(s)+(\dot\xi_1(s)+\dot\xi_2(s))^2\frac{\sin^2
2\eta(s)}{2}$ is constant along geodesics. Applying the initial condition $\eta(0)=0$, we
get
$$|v|^2=\dot\eta^2(0)=C-4\psi_1^2.$$ In the case when a geodesic
ends at the vertical line at $\xi_1(1)$,  its
lengths is expressed as
$$\sqrt{2}l_n=\sqrt{C-4\psi_1^2}=\sqrt{\pi^2n^2-4\xi_1^2(1)}.$$ We see that this result coincides with the result given
by Theorem~\ref{numbergeod1} and we state it as follows.

\begin{theorem}\label{numbergeod4}
Let $A$ be a point of the vertical line, i.~e. $\xi_1$ is given and  $\eta=0$. There are countably many geodesics $\gamma_n$ connecting $O=(0,0,0)$ and $A$, given parametrically as
\begin{eqnarray*}
\xi_1(s) &= & \xi_1s+\frac{1}{2}\arctan\Big(\frac{2\xi_1}{\pi n}\tan(\pi n s)\Big),\\
\xi_2(s) & = & -\xi_1 s,\\
\sin\eta(s) & = & \frac{s_n}{\pi n}\sin(\pi n s),
\end{eqnarray*} where $\frac{1}{\sqrt{2}}s_n=\frac{1}{\sqrt{2}}\sqrt{\pi^2 n^2-4\xi_1^2}$ is the length of the geodesic $\gamma_n$, $n\in\mathbb N$.
\end{theorem}

\section{Hopf fibration}

There is a close relation between the sub-Riemannian sphere $S^3$
and the Hopf fibration. Let $\mathbb S^2$ and $\mathbb S^3$ be unit 2-dimensional
and 3-dimensional sphere respectively. We remind that the Hopf
fibration is a principal circle bundle over two-sphere given by the map $h: \mathbb S^3\to \mathbb S^2$:
$$h(x_1,x_2,x_3,x_4)=((x_1^2+x_2^2)-(x_3^2+x_4^2),2(x_1x_4+x_2x_3),2(x_2x_4-x_1x_3)).$$
Another way to define the Hopf fibration is to write
$$h(q)=qiq^*\in \mathbb S^2,\quad q\in \mathbb S^3,\quad i=(0,1,0,0).$$ The fiber passing
through the unity of the group $(1,0,0,0)$ has equation
$(\cos\theta,\sin\theta,0,0)$, which as we see, coincides with the
equation of the vertical line at this point. The sphere $\mathbb
S^2$ represents the horizontal ``plane'' sweep out by the
geodesics with constant horizontal coordinates.
\begin{definition}
Let $Q\to M$ be a principle $G$-bundle with the horizontal distribution $\mathcal D$ on $Q$.
A sub-Riemannian metric on $Q$ that has distribution $\mathcal D$ and it is invariant under the action of $G$ is called a metric of bundle type.
\end{definition}

In our situation $\mathbb S^3\to \mathbb S^2$ is a principle $\mathbb S^1$-bundle given by the Hopf map.
The sub-Riemannian metric on the distribution
$\mathcal D=\spn\{X,Y\}$ was defined as the restriction of the euclidean metric $\langle\cdot,\cdot\rangle$
from $\mathbb R^4$ and we used the same notation $\langle\cdot,\cdot\rangle$ for sub-Riemannian metric.

\begin{proposition}
The sub-Riemannian metric $\langle\cdot,\cdot\rangle$ on $\mathbb S^3$ is a metric of bundle type.
\end{proposition}

\begin{proof}
The action of the group $\mathbb S^1$ on $q=(x_1,x_2,x_3,x_4)\in\mathbb S^3$ can be written as $q\circ e^{it}$,
$e^{it}=(\cos t+i\sin t)\in \mathbb S^1$, $t\in [0,2\pi)$, where $\circ$ is the quaternion multiplication.
If we write $q=(e^{i\xi_1}\cos\eta,e^{i\xi_2}\sin\eta)$, $\tilde q=q\circ e^{it}$, then
$\tilde q=(e^{i(\xi_1+t)}\cos\eta,e^{i(\xi_2-t)}\sin\eta)$. In order to show that the metric $\langle\cdot,\cdot\rangle$
is of bundle type, we have to prove that the metric is invariant under the action of the group $\mathbb S^1$.
The metric $\langle\cdot,\cdot\rangle$ at any $q$ is given by the matrix
\begin{equation*}
\left[\array{rrr}1 & 0 & 0
\\ 0 & \cos^2\eta & 0
\\
 0 & 0 & \sin^2\eta
\endarray\right]
\end{equation*} and it is easy to see that it is invariant under the action $\tilde q=q\circ  e^{it}$.
\end{proof}

We can formulate the results of Theorems~\ref{numbergeod1} and~\ref{numbergeod4} as an isoholonomic problem. First let us give some definitions. Let $c:[0,1]\to\mathbb S^2$ be a curve in $\mathbb S^2$ For a given point $x$ of the fiber $Q_{c(0)}$ at $c(0)$, let $\gamma$ be a horizontal lift of $c$ that starts at $q$ (it means that the projection of $\gamma$ under the Hopf map $h:\mathbb S^3\to\mathbb S^2$ coincides with $c$). The map $\Phi(c):Q_{c(0)}\to Q_{c(1)}$  sending $x=\gamma(0)$ to the endpoint $\gamma(1)$ of the horizontal lift is called {\it parallel transport} along $c$. The action of $\mathbb S^1$ takes horizontal curves to horizontal curves, so any two horizontal curves $\gamma_1$ and $\gamma_2$ of $c$ are related by $\gamma_1=\gamma_2g$ for some $g\in\mathbb S^1$. It follows that the action of $\mathbb S^1$ commutes with $\Phi(c)$, that is, $\Phi(c)(xg)=\big(\Phi(c)(x)\big)g$.

If $c$ is a closed loop, parallel transport $\Phi(c)$ maps the fiber $Q_{c(0)}$ onto itself. Fix a point $x_0\in Q_{c(0)}$. Since $\mathbb S^1$ acts transitively on $Q_{c(0)}$ we have $\phi(c)(x_0)=x_0l$ for some $l\in\mathbb S^1$. If we choose another point $y_0=x_0g\in Q_{c(0)}$, we get $$\Phi(c)(y_0)=\big(\Phi(c)(x_0)\big)(g)=x_0lg=y_0(g^{-1}lg).$$ The curve $c$ therefore, determines a conjugacy class in $\mathbb S^1$, called the holonomy class of $c$. The element $l\in\mathbb S^1$ for which $\Phi(c)(x_0)=x_0l$ is called the {\it representative holonomy} of $c$ with respect to $x_0$. The set of all such $l\in\mathcal S^1$ for $c$ running over all closed loops with $c(0)=c(1)=h(x_0)$ is a subgroup of $\mathbb S^1$ called the {\it holonomy group} of the distribution $\mathcal D$ at $x_0$.

Let us fix a representative holonomy $l\in\mathbb S^1$ and a point $x_0\in\mathbb S^3$, or equivalently the initial and the finite points $x_0$ and $x_1=x_0l$. The set of horizontal curves that join $x_0$ to $x_1$ is in one-to-one correspondence with the set of all closed loops on $\mathbb S^2$, based at $h(x_0)$, whose holonomy with respect to $x_0$ is $l$. Recall that the Riemannian length of a loop on $\mathbb S^2$ equals the sub-Riemannian length of its horizontal lift. Thus, the sub-Riemannian geodesic problem for geodesics with endpoints at the same fiber is equivalent to the following isoholonomic problem: {\it Among all loops with a given holonomy, find the shortest}. 

If we take $x_0=(1,0,0,0)$, $l=e^{i\omega}$, $x_1=x_0l=(\cos \omega,\sin \omega,0,0)$, then Theorem~\ref{numbergeod1} says 
\begin{eqnarray*}
\text{if}\quad \omega\in [0,\pi) & \quad\text{then the shortest loop has length}\quad s_1=\sqrt{\pi^2n^2-\omega^2},\\
\text{if}\quad \omega\in [\pi,2\pi) & \quad\text{then the shortest loop has length}\quad s_2=\sqrt{4\pi^2n^2-\omega^2}.
\end{eqnarray*}

\end{document}